\documentclass[12pt]{amsart}
\usepackage{amssymb}

\theoremstyle{definition}
\newtheorem{theorem}{Theorem}[section]
\newtheorem{lemma}[theorem]{Lemma}

\newtheorem{claim}[theorem]{Claim}

\newtheorem{definition}[theorem]{Definition}

\newtheorem{main lemma}[theorem]{Main Lemma}
\newtheorem{corollary}[theorem]{Corollary}

\newtheorem{question}[theorem]{Question}
\renewenvironment{proof}{\noindent\textbf{Proof}}
  {\hspace*{\fill} $\Box$ }

\newcommand{\rest}{\upharpoonright}
\newcommand{\into}{\rightarrow}
\newcommand{\embeds}{\hookrightarrow}

\newcommand{\QQ}{\mathbb{Q}}

\newcommand{\qkappa}{\QQ(\kappa)}

\newcommand{\ran}{\mbox{ran}}

\newcommand{\scf}{\mbox{\scriptsize{cf}}}

\newcount\skewfactor
\def\mathunderaccent#1#2 {\let\theaccent#1\skewfactor#2
\mathpalette\putaccentunder}
\def\putaccentunder#1#2{\oalign{$#1#2$\crcr\hidewidth
\vbox to.2ex{\hbox{$#1\skew\skewfactor\theaccent{}$}\vss}\hidewidth}}

\begin{document}

\title[The embedding structure for LOTS]{The embedding structure for linearly ordered topological spaces}
\author{A. Primavesi and K. Thompson}
\date{\today}

\begin{abstract}
In this paper, the class of all linearly ordered topological spaces (LOTS) quasi-ordered by the embeddability relation is investigated. In ZFC it is proved that for countable LOTS this quasi-order has both a maximal (universal) element and a finite basis. For the class of uncountable LOTS of cardinality $\kappa$ it is proved that this quasi-order has no maximal element for $\kappa$ at least the size of the continuum and that in fact the dominating number for such quasi-orders is maximal, i.e. $2^\kappa$. Certain subclasses of LOTS, such as the separable LOTS, are studied with respect to the top and internal structure of their respective embedding quasi-order. The basis problem for uncountable LOTS is also considered; assuming the Proper Forcing Axiom there is an eleven element basis for the class of uncountable LOTS and a six element basis for the class of dense uncountable LOTS in which all points have countable cofinality and coinitiality.
\end{abstract}

\maketitle

\section{Introduction}

The embedding structure for linear orders is dependent on the axiomatisation of set theory that we choose to adopt; this is a widely known example of set-theoretic independence phenomena occurring at the heart of classical mathematics. For example, there can consistently be a universal linear order (one that embeds every other linear order of the same cardinality) for every infinite cardinal $\kappa$. This statement can also consistently fail to hold. More recently, it has been shown that in certain models of set theory there is a five element basis for the uncountable linear orders, resolving a long-standing open question; other models of set theory do not admit a finite basis for the uncountable linear orders. 

A \emph{linearly ordered topological space} (or LOTS) is a linear order endowed with the open interval topology, 
call it $\tau$. 
An {\em embedding} is an injective function that preserves structure. 
Specifically, for the structures in this paper:

\begin{itemize}
\item A {\em linear
      order embedding} is an injective order-preserving map.
\item A {\em LOTS embedding} is an injective order-preserving map that is continuous.
\end{itemize}

The existence of a LOTS embedding, $f : A \rightarrow B$, where $A$ and $B$ are two arbitrary linear orders, ensures that not only is there a suborder of $B$, call it $B'$, that is order-isomorphic to $A$, but also that the open sets in $\tau_{B'}$ are exactly those sets of the form $B' \cap u$ for some $u \in \tau_B$. The existence of a LOTS embedding is a non-trivial, natural strengthening of the existence of a linear order 
embedding.

The relation on the class of all LOTS defined by setting $A \preceq B$ if and only if $A$ LOTS-embeds into $B$ is a quasi-order. We note that it is not a partial order because bi-embeddability (where $A \preceq B$ and $B \preceq A$ both hold) does not imply isomorphism.  Henceforth, when we talk about the `embeddability order' we will generally be referring to the quasi-order $\preceq$. However, when we discuss e.g.\ chains in this order, we will be referring to the strict ordering $\prec$. 
As we will see, the quasi-order of LOTS embeddability in general looks very different from that of linear order embeddability, but by restricting our attention to certain subclasses of LOTS we can obtain similarities.

In this paper we investigate the properties of the embedding structure for LOTS and compare them to the properties of the embedding structure for linear orders. 

There are several aspects of the LOTS embeddability order that 
we study: here we split them into three groups. (1) The top (the question of \emph{universality}); can we find a non-trivial set of linear orders such that every linear order in a given class must be LOTS-embeddable into one of them? (2) The bottom (the \emph{basis} question); can we find a non-trivial set of linear orders such that every linear order in a given class must LOTS-embed one of them? (3) The internal structure of the embeddability order; what are the possible cardinalities of chains and antichains 
in the LOTS embeddability order for a given class of LOTS?

We include results pertaining to all three questions. The structure of this paper is as follows:

In Section \ref{sec:prelims} we review some basic notation and explain the various models of set theory in which we will be proving our results.

In Section \ref{sec:internal} we study the internal structure of the LOTS embedding order, proving that it has large chains and antichains when the LOTS have cardinality at least the continuum.
In Section \ref{sec:universal} we investigate the {\em top} of the embeddability quasi-order. We show that the rationals form a universal for countable LOTS, but that there can be no universal LOTS for any cardinal $\kappa$ greater than or equal to the size of the continuum; in fact, we prove that the dominating number for the LOTS embeddability order is maximal, i.e. $2^\kappa$.
We can however define a class of linear orders of size $\kappa$ for which a universal LOTS exists, called the $\kappa$-entwined LOTS. 

By restricting to certain subclasses of LOTS we can find similarities between the LOTS embeddability order and that of the linear orders.
In Section \ref{sec:dense} we prove that universal LOTS for the class of separable partial orders can exist at various cardinals, by way of some general results on dense subsets and continuity. These results generalise to other cardinals, with an appropriate generalisation of the notion of separability.
An analogue of a theorem of Sierpinski shows that there is a sequence of length continuum of separable LOTS that are strictly 
decreasing in the LOTS embeddability order.

In Section \ref{sec:basis} we investigate the {\em bottom} of the embeddability order. Our results in this section are obtained under the assumption of PFA, the Proper Forcing Axiom. It has recently been proved by J. Moore \cite{5-elements} that in
models of set theory where the PFA holds, there is a five-element
basis for the uncountable linear orders. We prove that under PFA there is an 11 element basis for uncountable LOTS, and show that this is the smallest possible. We also prove that there is a six element LOTS basis for linear orders that are dense and have only points of countable cofinality/coinitiality.

We conclude in Section \ref{sec:questions} with some open questions.


\section{Notation and preliminaries}\label{sec:prelims}

In this paper we work under the standard axiomatic assumptions of ZFC (see e.g.\ \cite{jech} for details) and certain additional axioms. That is, we assume the consistency of ZFC and work within a model of set theory in which the sets obey these axioms. When the axioms that are being assumed for a given proof are not explicitly stated, it is to be understood that the proof requires only the axioms of ZFC. However as the additional axioms that we use are together inconsistent (but both consistent with ZFC), our results show that ZFC is often not enough. 


The first additional axiom that we use is the Continuum Hypothesis (CH) and its generalisation GCH, which states that for any infinite cardinal $\kappa$, $2^\kappa = \kappa^+$. That is, the powerset of $\kappa$ has the smallest possible cardinality. In this notation CH is written as $2^{\aleph_0} = \aleph_1$.

The second is the Proper Forcing Axiom (PFA). This is a strong Baire category assumption that can be proved consistent via forcing, assuming the consistency of ZFC and certain large cardinal axioms. PFA is extremely useful for proving independence results concerning linear order embeddings, as we will see in this paper. In models where PFA holds it is provable that $2^{\aleph_0} = \aleph_2$, so in particular CH fails. For more on this see \cite{Moore:PFAsurvey}.

We also use methods from standard topology, especially the topology of the reals. 
For general topological definitions see e.g.\ \cite{munkres}. We will now introduce some of the notation used throughout the paper:

If $a$ is an element and $A \subseteq L$ a subset of a linear order, then we write $a < A$ when $a < b$ 
for all $b \in A$. For $A, B \subseteq L$, $(A, B)$ denotes the convex set $\{ x \in L: \forall y \in A (y <_L x) \text{ and } \forall y \in B (x <_L y)\}$. If $f$ is a function with domain $X$ and $x \subseteq X$ then we denote by $f[x]$ the pointwise image of $f$ under $x$.

For $x \in L$, let $l(x) =  \{ y \in L : y < x \}$ and $r(x) = \{ y \in L : x < y \}$. The {\em cofinality} of $x$ is the cardinality of the smallest increasing sequence for which $x$ forms a supremum, and the {\em coinitiality} is the cardinality of the smallest decreasing sequence for which $x$ forms an infimum; of course, the cofinality/coinitiality function is not necessarily defined for all $x \in L$. For two linear orders $A$ and $B$, $A \times B$ will always denote the lexicographical product of $A$ with $B$, and $A^*$ denotes the reverse ordering of $A$, $(A, >_A)$.

Formally, a LOTS is denoted by $(L, <_L, \tau_L)$ where $L$ is a set with a linear ordering $<_L$ and $\tau_L$ is the open interval topology on the linear order $(L, <_L)$. However, we usually drop the subscript from the linear ordering when it is clear from context and we use $\tau_L$ implicitly. In fact, we will usually denote such LOTS simply by $L$.

Strictly speaking `LOTS' is a singular term; following standard usage we will use it for the plural as well. 

\section{The internal structure of the LOTS embeddability order}\label{sec:internal}

The internal structure of the LOTS embeddability order was considered by Beckmann, Goldstern and
Preining in \cite{cont-fraisse} where they show that countable closed
sets of reals are well-quasi-ordered by LOTS embeddability. That is, 
they show that the embedding quasi-order for this class of LOTS has
no infinite decreasing sequence and no infinite antichain. 

In this section we show that uncountable LOTS are not well-behaved with respect to their embedding structures. In particular, there are large antichains and chains. We adapt methods from standard topology of the reals, namely linear continua and the intermediate value theorem. 

The intermediate value theorem (I.V.T.) tells us that a continuous embedding from a linear continuum, $A$, into a linear order $B$ must be surjective onto a convex subset of $B$, and therefore that there is a convex subset of $B$ isomorphic to $A$. In Section \ref{sec:universal} we will make use of this fact to prove that the top of the embedding quasi-order is maximally complex for LOTS of size $\kappa$, where $\kappa \geq 2^\omega$ (because no linear continuum can exist with cardinality less than $2^\omega$).

It is easy to define several linear continua of size equal to the continuum. We use them to `code' subsets of $\kappa$, as follows:

\begin{lemma}

Let $[0, 1] \subseteq \mathbb{R}$ denote the closed unit interval -- that is, a copy of $\mathbb{R}$ with endpoints -- and $[0, 1)$ an isomorphic copy of $\mathbb{R}$ with a least point but no greatest point. Then each of the following is a linear continuum:

\begin{itemize}

\item $[0, 1]$.

\item $[0, 1)$.

\item $R^\prime = [0, 1) \times [0, 1]$.

\item $R^{\prime\prime} = [0, 1) \times [0, 1] \times [0, 1]$.

\item $R_0 = R^\prime + [0, 1)$.

\item $R_1 = R^{\prime\prime} + [0, 1)$.

\end{itemize}

\end{lemma}

\begin{proof}
The first two are trivial. For the rest, recall that when dealing with linear orders of the form $I \times J$, ordered lexicographically, any infinite increasing (or decreasing) sequence in this ordering will have an infinite subsequence either entirely contained within $\{ i \} \times J$, for some $i \in I$, or with the property that for any $i \in I$ there is at most one $j \in J$ with $(i, j)$ appearing in this subsequence. In the former case, the existence of a supremum (or infimum) to this subsequence then follows if $J$ is a linear continuum, and $J$ has endpoints. In the latter case it follows from the existence of a supremum (or infimum) in $I$ to the set of $i \in I$ such that $(i, j)$ is in this subsequence for some $j \in J$, and the fact that $J$ has a greatest and least point. A modification of this argument gives us the required proof in each case.
\end{proof}

We will use an infinite sum of copies of $R_0$ and $R_1$ to code subsets of $\kappa$. The following is apparent:

\begin{lemma}

$R_0$ cannot be continuously embedded into $R_1$. Likewise, $R_1$ cannot be continuously embedded into $R_0$.
\end{lemma}

\begin{proof}
By the I.V.T., if there was such an embedding then $R_1$ would contain an interval isomorphic to $R_0$, or vice versa. This is clearly not the case.
\end{proof}

We also remark that because $R_0$ and $R_1$ both have a least point, any direct sum of the form $\sum_{\alpha < \zeta} R_{i_\alpha}$, where $i_\alpha \in \{ 0, 1 \}$ and $\zeta$ is an ordinal, is also a linear continuum. We are now ready to prove the following:

\begin{theorem}\label{thm:chains}

Let $\kappa \geq 2^{\aleph_0}$. Then there exists:

\begin{enumerate} 

\item A set $\mathcal{A}$ of LOTS of size $\kappa$, with $| \mathcal{A} | = 2^\kappa$, such that there is no LOTS embedding from $A$ into $B$ for any two distinct $A, B \in \mathcal{A}$.\\

\item A sequence $I = \langle I_\zeta : \zeta < \kappa^+ \rangle$ of LOTS of size $\kappa$, strictly increasing in the LOTS embeddability order.\\

\item For each $\eta < \kappa^+$, a sequence $D^\eta = \langle D_\zeta : \zeta < \eta\rangle$ of LOTS of size $\kappa$, strictly decreasing in the LOTS embeddability order.

\end{enumerate}

\end{theorem}

\begin{proof} 
Let $X$ and $Y$ be subsets of $\kappa^+$ of size $\kappa$. Define $\alpha_X$ to be $\mathrm{sup}\{ \beta + 1 : \beta \in X \}$ and similarly define $\alpha_Y$. Let $g_X : \alpha_X \rightarrow 2$ be the characteristic function of $X$, and $g_Y$ be the same for $Y$. As $X$ and $Y$ are bounded in $\kappa^+$ we must have $\alpha_X, \alpha_Y < \kappa^+$. Let $R_X = \sum_{\alpha < \alpha_X} R_{g_{X}(\alpha)}$ and $R_Y = \sum_{\alpha < \alpha_Y} R_{g_{Y}(\alpha)}$. Both $R_X$ and $R_Y$ are linear continua for any such $X$ and $Y$. And because $\kappa \geq 2^\omega$, both are of size $\kappa$. 

By the I.V.T. if there were a LOTS embedding from $R_X$ to $R_Y$ then there would have to be a convex subset of $R_Y$ order isomorphic to $R_X$. Whenever $X \neq Y$ and $\alpha_X = \alpha_Y = \kappa$, this will only happen if the values of the characteristic function of $X$ (considered as an uncountable string of 0's and 1's) are equal to a final section of the values of the characteristic function of $Y$. But if $X$ and $Y$ are such that this is not the case, then the I.V.T. tells us that $R_X$ and $R_Y$ are pairwise non-embeddable as LOTS.

(1) We can find a family of $2^\kappa$ many subsets of $\kappa$, $\{ X_\zeta : \zeta < 2^\kappa \}$, such that each is unbounded in $\kappa$ and for any two, $X_\beta$, $X_\gamma$, the characteristic function of $X_\beta$ (as above) is not equal to a final section of the characteristic function of $X_\gamma$; we can always find such a family by diagonalisation. Thus $\mathcal{A} = \{R_{X_\zeta} :  \zeta < 2^\kappa\}$ is an antichain in the LOTS embedding order.

(2) Let $\langle Y_\zeta : \zeta < \kappa^+ \rangle$ be a sequence of bounded subsets of $\kappa^+$ such that for each $\zeta < \zeta' < \kappa^+$ we have $Y_\zeta$ is a proper initial segment of $Y_{\zeta'}$. Then let $I_\zeta = R_{Y_\zeta}$ for $\zeta < \kappa^+$. 

(3) Let $\eta < \kappa^+$ and $Y$ be an unbounded subset of $\eta$ such that there is a set $Y' \subseteq \eta$ with any two distinct $\alpha, \beta \in Y'$ being such that $\langle g_Y(i) : \alpha < i < \eta \rangle \neq \langle g_Y(i) : \beta < i < \eta \rangle$ and $\mathrm{otp}(Y') = \eta$. Then let $R_\zeta = \sum_{\zeta < \alpha < \eta} R_{g_{Y}(\alpha)}$, for each $\zeta \in Y'$. Then $D^\eta = \{ D_\zeta = R_\zeta : \zeta \in Y' \}$ is as required.
\end{proof}

\section{Universality for general LOTS}\label{sec:universal}

For a given class of structures a {\em universal} at cardinality $\kappa$ is a member of the class, of size $\kappa$,
that embeds all the other structures in this class having size $\kappa$. Universals were first studied in
topology, where the investigation of these quintessential objects was used to determine properties of the class of structures as a whole. The study was broadened using model theory, which tells us that under GCH every first order theory in a countable language has a universal in all uncountable cardinals. The class of linear orders has this property. However, the model theoretic technique that determines this cannot say anything regarding LOTS, as the theory of these structures is not definable in first-order logic. In this case we must use set-theoretic techniques to determine questions of universality. 

For regular cardinals, the universal linear orders that exist under GCH are natural and can be easily constructed. 
At cardinality $\aleph_0$, $\QQ$ is the universal linear order. For an uncountable regular
cardinal, $\kappa$, there is a universal linear order $\QQ(\kappa)$
that has a generalisation of the density property of the rationals,
called $\kappa$-saturation:

\[ \forall S,T\in [\QQ(\kappa)]^{<\kappa} \hspace{.1cm} [S <
	T\implies (\exists x) S< x < T].\]

As with the rationals, $\QQ(\kappa)$ is the unique (up to isomorphism)
$\kappa$-saturated linear order of size $\kappa$ without endpoints. It
exists whenever $\kappa = \kappa^{< \kappa}$, so in particular under GCH. It
was known to exist by Hausdorff in 1908 (\cite{haus}) and
was constructed explicitly by Sierpinski using lexicographically ordered
sequences of 0's and 1's of length $\kappa$ which have a final 1 (see \cite{linorders} for
details). If $\kappa$ is a singular cardinal and $\kappa = \kappa^{< \scf(\kappa)}$ then there is also
a universal linear order which is not saturated but is still \emph{special} (it is the direct limit of saturated
linear orders of smaller (regular) cardinality, see e.g.\ \cite{kies:chang} for details). 

We note that $\kappa$-saturation for regular uncountable $\kappa$ tells us that no infinite increasing sequence of cofinality less than $\kappa$ can have a supremum
in $\QQ(\kappa)$, so it cannot possibly be a universal for LOTS of size $\kappa$. Specifically, any ordinal of the form $\alpha + 1$ for $\alpha$ a limit ordinal less than $\kappa$ cannot continuously embed into it (and clearly we can find a LOTS of size $\kappa$ with an interval isomorphic to $\alpha + 1$). Taking the completion of $\QQ(\kappa)$ under sequences of length less than $\kappa$ will negate such counterexamples and thus gives a universal for a broad class of LOTS called the $\kappa$-{\em entwined} LOTS (see definition below), but it is still not a universal for general LOTS because in this case $\alpha + 1 + \alpha^*$ cannot continuously embed into it, where $\alpha$ is as above. Similar counterexamples show that the full Dedekind completion of $\qkappa$ is not a universal for LOTS of size $2^\kappa$ and likewise for the special linear orders that exist under GCH (particularly at singular $\kappa$). 

\begin{definition}

A LOTS, $L$, of size $\kappa$ is $\kappa$-{\em entwined} if for all $x \in L$, sup$(l(x))$ = inf$(r(x))$ = $x$ implies that both the cofinality and coinitiality of $x$ are equal to $\kappa$.

\end{definition}

\begin{theorem}
\label{entwined}
(GCH)
Let $\bar{\QQ}(\kappa)$ be the completion of $\QQ(\kappa)$ under sequences of length $<\kappa$. Then $\bar{\QQ}(\kappa)$ has size $\kappa$ and is universal for $\kappa$-entwined LOTS.
\end{theorem}

\begin{proof} Using the linear order universality of $\QQ(\kappa)$, we
  construct the relevant embeddings by taking linear order embeddings and
  altering them so as to be
  continuous. Let $L$ be a $\kappa$-entwined linear order of size $\kappa$ and let $f :
  L \embeds \QQ(\kappa)$ be an injective order-preserving map. 

Call a point $x \in L$ a $\kappa$-point if it is the infimum or supremum of a decreasing/increasing sequence of length $\kappa$ in $L$.

We say that the function $f$ is {\em discontinuous} at a point $x \in L$ if and only if
$x$ = sup($l(x))$ but $f(x) \neq \text{sup}(f[l(x)])$ or $x$ = inf($r(x))$ but $f(x) \neq \text{inf}(f[r(x)])$. We start by ``fixing'' the discontinuity of $f$ for those points in $L$ that are $\kappa$-points.

So, for every $\kappa$-point
$x$ at which $f$ is discontinuous (denote the set of such points by $D$), we remove all the members of $\QQ(\kappa)$ contained in the convex set between $f(x)$ and $f[l(x)]$, and between $f(x)$ and $f[r(x)]$.

We will call the thinned order $\QQ'(\kappa)$. 
Let 
\[\QQ'(\kappa) = \QQ(\kappa) \setminus \bigcup_{x \in D} \{y \in \QQ(\kappa) : y \in (f[l(x)], f(x)) \cup (f(x), f[r(x)]).\]
We shall show that this is isomorphic
to $\QQ(\kappa)$ and thus the composition
of this isomorphism (which is trivially continuous) with $f$ will not be discontinuous at any $\kappa$-point in $L$.

Firstly, observe that 
only convex sets between points in the range of
$f$ were removed from $\QQ(\kappa)$, so $\QQ'(\kappa)$ has no endpoints. It is also clear that $\QQ'(\kappa)$ has size $\kappa$, as it contains an isomorphic copy of $L$.

We will show that $\QQ'(\kappa)$ is $\kappa$-saturated. Let $A, B \subseteq \QQ'(\kappa)$ be such that both are sets of size less than $\kappa$ and $A <
B$ and assume towards a contradiction that $(A, B)_{\QQ'(\kappa)} = \emptyset$. In
this case, we must have that $\ran(f) \cap (A, B)_{\QQ(\kappa)} = \emptyset$
because $\QQ'(\kappa)$ contains the entire range of $f$.

The convex set $(A, B) \subseteq \QQ(\kappa)$ must be of the form $(f[l(x)], f(x))$ or $(f(x), f[r(x)])$ for some $x \in D$, by the construction of $\QQ'(\kappa)$. But in either case one of $A$ or $B$ must have cofinality or coinitiality $\kappa$ because $x$ is a $\kappa$-point, contradicting our assumption that both $A$ and $B$ have size less than $\kappa$.

So $L$ is embedded into $\QQ'(\kappa)$ by a map, $f$, that is continuous at $\kappa$-points; $\QQ'(\kappa)$ is a $\kappa$-saturated linear order of size $\kappa$ without endpoints, and hence is isomorphic to $\QQ(\kappa)$. This establishes the existence of a map $f' : L \rightarrow \QQ(\kappa)$ that is continuous at $\kappa$-points.
We define a function $g: L \rightarrow \bar{\QQ}(\kappa)$ that is continuous everywhere. Let $f'$ be discontinuous at a point $x \in L$, then $x = \text{sup}(l(x))$ or $x = \text{inf}(r(x))$, but not both because $L$ is $\kappa$-entwined. $\bar{\QQ}(\kappa)$ is complete for sequences of length $< \kappa$, so in the case that $x = \text{sup}(l(x))$ we set $g(x) = \text{sup}(f'[l(x)])$ and if $x = \text{inf}(r(x))$ we set $g(x) = \text{inf}(f'[r(x)])$, which is possible in both cases because $x$ is not a $\kappa$-point. At all other points $x$ we let $g(x) = f'(x)$.

Our construction of $g$ ensures that it is both order-preserving and continuous. Hence $\bar{\QQ}(\kappa)$ is a universal for $\kappa$-entwined LOTS.

\end{proof}

\begin{corollary}
\label{Qisuniversal}
Every countable LOTS is $\omega$-entwined, so the rationals $\QQ$ is a universal for countable LOTS.
\end{corollary}

Corollary \ref{Qisuniversal} is true in ZFC, but Theorem \ref{entwined} only applies in general when we assume GCH. Without GCH, the situation is more complicated. 
The definition of a LOTS embedding tells us that if there
is no universal linear order at a particular cardinal, then there
cannot be a universal LOTS at that cardinal. So we will only search for universal LOTS in
those cardinalities, models of set theory and subclasses where 
universal linear orders are known to exist. 
In the absence of GCH we are constrained by the following Theorem:

\begin{theorem}[Kojman-Shelah \cite{shel:koj}]
\label{kojmanshelah}
For any regular $\kappa \in (\aleph_1, 2^{\aleph_0})$ there is
no universal linear order of size $\kappa$. 

\end{theorem}

Kojman and Shelah's proof extends to show that club 
guessing at $\aleph_1$ (which is a weak combinatorial principle consistent 
with ZFC, see e.g.\ \cite{mirna-survey}) and a failure of CH 
imply that there is no universal linear order at $\aleph_1$. This fact, together with Theorem \ref{thm:no-univ} below establishes the non-existence of universal LOTS at almost all cardinalities; the countable case is a rare exception.

\begin{theorem}\label{thm:no-univ} There cannot be a universal LOTS at any cardinality $\kappa \geq 2^{\aleph_0}$.
\end{theorem}

In fact, we will prove something stronger, which intuitively says that the top of the LOTS embedding order is maximally complex: 

\begin{definition}

The {\em dominating number} for a quasi-order $P$ is the least possible cardinality of a subset $Q \subseteq P$ such that for any $p \in P$ there is a $q \in Q$ with $p \leq q$. 

\end{definition}

If a universal linear order or LOTS exists then the appropriate embeddability order has a dominating number of 1. We will show that the dominating number for the embeddability order for LOTS of size $\kappa$ is the maximum possible, namely $2^{\kappa}$, whenever $\kappa \geq 2^\omega$.

\begin{theorem}
\label{noLOTS}
Let $\kappa \geq 2^\omega$. Then the dominating number for the embeddability order for LOTS of size $\kappa$ must be $2^{\kappa}$. Consequently, there can be no universal LOTS of size $\kappa$. 

\end{theorem}

\begin{proof}
We prove the first statement; the second is an immediate consequence of it. Assume towards a contradiction that $\mathcal{U}$ is a family of size $\lambda < 2^\kappa$ and witnesses the fact that the dominating number for the LOTS embeddability order is $\lambda$. Let $\mathcal{A}$ be the antichain asserted to exist in Theorem \ref{thm:chains} (1). $\mathcal{A}$ has cardinality $2^\kappa$, hence we can find a $U \in \mathcal{U}$ such that $2^\kappa$ many elements of $\mathcal{A}$ continuously embed into $U$. By the construction of $\mathcal{A}$, every member of it is a linear continuum so all those members of $\mathcal{A}$ that LOTS embed into $U$ must be isomorphic to a convex subset of $U$. But $\mathcal{A}$ is an antichain under LOTS embeddability, so any two such convex subsets must be disjoint. Hence $U$ contains $2^\kappa$ many non-empty disjoint convex subsets, contradicting its size being $\kappa$.
\end{proof}

In light of Theorems  \ref{kojmanshelah} and \ref{noLOTS}, we now have the following:

\begin{corollary}

If $\kappa$ is an uncountable cardinal and a universal LOTS of size $\kappa$ exists, then $\kappa < 2^\omega$ and $\kappa = \aleph_1$ or is singular. 

\end{corollary}

Shelah, in \cite{sh100}, has proved it consistent that there is a universal linear order at $\aleph_1$ for the case where $\aleph_1 < 2^{\aleph_0}$. However, this linear order is not a universal LOTS. We do not know if there can consistently exist a universal LOTS at any uncountable cardinal. In general, models of a failure of CH do not have universal linear orders of size $\aleph_1$. In particular,  Shelah notes (see \cite{shel:koj}) that adding
$\aleph_2$-many Cohen reals over a model of CH produces a model in
which there does not exist a universal linear order of size
$\aleph_1$.

The proof of Theorem \ref{noLOTS} makes essential use of the fact that a linear continuum is topologically {\em connected}. It therefore establishes that there can be no universal of cardinality $\kappa \geq 2^\omega$ for the subclass of connected LOTS. A natural question then is: can we get the same result for the class of disconnected LOTS? We can easily construct disconnected spaces to play the role of the linear continua in the previous proof, by adding a pair of isolated points to each LOTS constructed in Theorem \ref{thm:chains}(1), so instead we ask a stronger question:

\begin{definition}
A linear order $L$ is {\em densely disconnected} if for any distinct $a, b \in L$ there is a partition of $L$ into two disjoint open sets, both of which have non empty intersection with $[a, b]$.

\end{definition}

The definition of a densely disconnected linear order is related to the usual topological notion of a totally disconnected space.

\begin{theorem}\label{thm:noLOTSc} Let $\kappa$ be a regular infinite cardinal. Then there is no universal for the subclass of densely disconnected LOTS of size $2^{\kappa}$, assuming that 
$\kappa^{< \kappa} = \kappa$.
\end{theorem}

\begin{proof} Let $\mu = 2^{\kappa}$. For a contradiction, assume that $L = (\mu, <_L, \tau_L)$ 
is a universal for the densely disconnected LOTS of size $\mu$. We will construct $A$, a densely disconnected LOTS with underlying set $\mu$, by induction, such that $A$ does not continuously embed into $L$.

Let $A \rest \kappa \cong \qkappa \times 2$. 

At stage $\alpha < \mu$ in the induction, assume that we have constructed $<_{A \rest c_\alpha}$ for some 
ordinal $c_\alpha \in [\kappa, \mu)$ but that $<_{A \rest c_\alpha +1}$ is undefined. 

Each stage of this inductive construction will take care of an embedding from 
$A \rest \kappa$ into $L$, so that all such possible embeddings cannot be extended to one that is both order-preserving and continuous from $A$ into $L$.

To this end, let $f_\alpha : A\rest\kappa \into L$ be an injective order-preserving function which has not yet been considered. (We may assume a canonical ordering of these functions so that we are assured to hit all of them in an induction of order type $\mu$.) Choose $I_\alpha = \langle i^\alpha: \alpha < \kappa\rangle , D_\alpha = \langle d^\alpha : \alpha < \kappa\rangle$ to be bounded sequences of $A \rest \kappa$ such that
\begin{enumerate} \item for all $\beta < \kappa$ we have $i^\alpha_\beta <_A i^\alpha_{\beta+1}$ and $d^\alpha_\beta >_A d^\alpha_{\beta+1}$
\item $I_\alpha <_A D_\alpha$
\item there does not exist an $x \in c_\alpha$ such that $I_\alpha <_{A\rest c_\alpha} x  <_{A\rest c_\alpha} D_\alpha$.
\end{enumerate}

$\QQ(\kappa)$ is $\kappa$-saturated so we can always choose such sequences because otherwise our order $A \rest c_\alpha$ would have to be complete and thus would already have size $\mu$, but $c_\alpha < \mu$ so this is not possible.

We will extend the ordering on $A$ as follows: 
If there is an $x \in L$ such that sup$(f_\alpha[I_\alpha]) = x$ = inf$(f_\alpha[D_\alpha])$ then let $I_\alpha <_A c_\alpha <_A c_{\alpha+1} <_A D_\alpha$ and take the transitive closure. If there is no such $x$, then we extend the ordering so that $I_\alpha <_A c_\alpha <_A D_\alpha$ and again take the transitive closure. 

If $A \rest c_\beta$ is defined for all $\beta < \alpha$ limit, then let $<_{A \rest c_\alpha} = \bigcup_{\beta < \alpha} <_{A \rest c_\beta}$ and let $<_A =  \bigcup_{\beta < \mu} <_{A \rest c_\beta}$.

Now let $f : A \into L$ be a LOTS embedding, which exists because $A$ has size $\mu$ and we are assuming that $L$ is universal at $\mu$. Then the initial segment $f \rest \kappa$ was considered in the inductive construction as $f_\delta$ for some $\delta < \mu$. At this stage of the induction, we chose sequences $I_\delta, D_\delta$. If sup$(f_\delta[I_\delta]) = x = $inf$(f_\delta[D_\delta])$ for some $x$ in $L$, then we set two distinct points in between $I_\delta$ and $D_\delta$ which clearly contradicts either the continuity or order-preservation properties of $f$.

The remaining cases to consider are that when sup$(f_\delta[I_\delta]) \neq$ inf$(f_\delta[D_\delta])$ both exist in $L$ and those where either sup$(f_\delta[I_\delta])$ or inf$(f_\delta[D_\delta])$ does not exist in $L$. In this case we added a unique point $x$ to $A$ such that sup($I_\delta) = x = \text{inf}(D_\delta)$. Again, this clearly contradicts the possibility that $f$ could be both order-preserving and continuous. Sequences $I_\gamma$ and $D_\gamma$ at later stages of the induction were chosen so that we would not add any further points between $I_\delta$ and $D_\delta$, by item (3) in the definition of $I_\gamma, D_\gamma$. This shows that $A$ cannot be continuously embedded into $L$ and therefore contradicts the universality of $L$, so long as $A$ is densely disconnected.

To see that this is indeed the case, let $a, b \in A$ be distinct and assume that $[a, b]$ cannot be split into two disjoint open intervals. Then $A\rest\kappa \cap [a,b]$ must be empty, but this is clearly not possible as only single points or pairs of points were added between elements of $A\rest\kappa$; in the latter case $[a, b]$ can be split into two disjoint open intervals and the former case contradicts our assumptions on $a,b$.

\end{proof}





\section{Dense embeddings and separable LOTS}\label{sec:dense}

In this section we note that a special kind of linear order embedding implies continuity and use this to show similarities between the embedding structures for certain subclasses of linear orders and LOTS. In particular, embeddings between separable linear orders, when they exist, have this property. Thus, we show that the embedding structure of separable LOTS, unlike in the general case, bears a close resemblance to the linear order embedding structure.

\begin{definition}\label{dense-emb}

Let $A$ and $B$ be linear orders and $f : A \rightarrow B$ a linear order embedding. Then we call $f$ a {\em dense embedding} if there is a convex subset $C \subseteq B$ such that $f[A]$ is a dense subset of $C$.

\end{definition}

If there exists a universal LOTS, $U$, for a subclass of linear orders, then 
taking the Dedekind completion of $U$ will give a universal LOTS
for those linear orders which densely embed $U$. Of course, the completion will not in general be of the same cardinality as $U$. The following simple lemma shows that with dense order-preserving embeddings we get continuity for free:

\begin{lemma}\label{dense-cont} Let $A, B$ be linear orders and suppose 
$f : A \into B$ is a dense embedding.
Then $f$ is continuous. 
\end{lemma}

This follows easily from the fact that the existence of any point of discontinuity $x$ implies there are elements of $C$ (as in Definition \ref{dense-emb}) between the image of a sequence and its limit $f(x)$, and density implies that some of these are in $f[A]$ contradicting the order preservation of $f$.

\begin{theorem}  Let $\mathcal{A}$ be a class of linear orders
each of size $\kappa$ and suppose that  $\mathcal{A}$ has a universal
LOTS. Then there is a universal LOTS for $\mathcal{A}'$, the class of all linear orders $L'$ such 
that there is an $L \in \mathcal{A}$ which densely embeds into $L'$.
\end{theorem}

\begin{proof} Let $U$ be universal for $\mathcal{A}$ and take the completion of 
$U$ under sequences of length $\leq \kappa$, call this $\bar{U}$. Note that
$U$ is dense in $\bar{U}$. 

For any linear ordering $L' \in \mathcal{A}'$ there is an $L \in \mathcal{A}$
that is dense in $L'$. One may then
find an order-preserving continuous $f : L \into U$ and a dense embedding $f' : U \into \bar{U}$. The composition of these two functions gives a continuous map from $L$ into $\bar{U}$. Since $L$ is a dense subset of $L'$ this induces a continuous map from $L'$ into $\bar{U}$ by the definition of $\bar{U}$. 
\end{proof}

\begin{corollary}\label{dense-examples} There is a universal for all separable LOTS, namely 
the reals. Moreover, if $\kappa = \kappa^{< \kappa}$ then there is a 
universal LOTS for linear orders of size $2^\kappa$ which densely embed orders of size $\kappa$
that are $\kappa$-entwined.
\end{corollary}

The existence of dense embeddings for separable linear orders produces very similar embeddability results for LOTS as for linear order embeddings, as in the case of universality. We contrast what happens at the cardinality of the continuum with the situation for separable LOTS of cardinality $\aleph_1$ under PFA, where in particular $\aleph_1 < 2^{\aleph_0}$. 

The next results follow from a theorem of Sierpinski, see \cite[Theorem 9.10]{linorders}:

\begin{theorem} $\text{}$
\begin{itemize} 
\item There is a sequence $\langle X_\alpha : \alpha < 2^{\aleph_0}\rangle$ of 
separable linear orders each of size $2^{\aleph_0}$ which is strictly decreasing in the LOTS 
embedding quasi-order. That is, for every $\beta < \alpha  < 2^{\aleph_0}$, $X_\alpha$ 
LOTS embeds into $X_\beta$ but $X_\beta$ does not LOTS embed into $X_\alpha$.
\item There is a set $\{ Y_\alpha : \alpha < 2^{2^{\aleph_0}}\}$ of 
separable linear orders each of size $2^{\aleph_0}$ which are pairwise incomparable in the embeddability order for linear orders and therefore also for LOTS. 
\end{itemize}
\end{theorem}

The decreasing sequence and antichain that Sierpinski constructed consist of dense sets of reals. The
embedding from $X_\alpha$ into $X_\beta$ for $\beta < \alpha$ is always the identity so is trivially continuous. The fact that there is no linear order embedding between any $Y_\alpha, Y_\beta$ for $\alpha \neq \beta$ implies that there is also no LOTS embedding.

 A result of Baumgartner in \cite{baum} shows that in models of PFA all $\aleph_1$-dense sets of reals
are isomorphic. This set of reals has size $\aleph_1$, no endpoints and has the property that
any interval has size $\aleph_1$. Every separable LOTS of size $\aleph_1$ is isomorphic to a set of reals (by the universality of the rationals) and can be extended to an $\aleph_1$-dense set of reals.
This gives us the following:

\begin{theorem}
\label{pfa-separable} Under PFA, there is a universal separable LOTS of size $\aleph_1$. This is also a basis for such LOTS.
\end{theorem}

\begin{proof} There is a continuous embedding from the countable dense subset of a given separable LOTS to the rationals, and this induces a continuous map from the LOTS itself to a canonically chosen $\aleph_1$-dense set of reals, by Baumgartner's result.

Every separable linear order of size $\aleph_1$ has an $\aleph_1$-dense suborder. Thus the unique $\aleph_1$ dense set of reals (up to isomorphism) embeds densely into this convex set. So there is a LOTS embedding from the $\aleph_1$ dense set of reals into any separable linear order of size $\aleph_1$. 

\end{proof}

Under PFA, the $\aleph_1$-dense set of reals forms a single element 
basis for the uncountable separable linear orders and also a universal for those of size $\aleph_1$. 
This makes the embedding quasi-order for separable linear orders of size $\aleph_1$ 
completely flat since they are all bi-embeddable. By Theorem \ref{pfa-separable} the same is true for LOTS.

In summary, the separable linear orders and LOTS of size $\aleph_1$ have a universal assuming CH but have a rather chaotic internal embedding structure involving long $\prec$-chains and antichains with as many orders as there are isomorphism classes. However, in certain situations without CH (namely when PFA holds), the embedding structure consists of a single embeddability class.

\section{The basis question for LOTS}\label{sec:basis}

The countably infinite linear orders have a two element basis, consisting of $\omega$ and $\omega^*$. This is trivially also a basis for the countable LOTS. 

Under CH there can be no finite (or even countable) basis for the uncountable linear orders, by a result of Dushnik and Miller (\cite{Dushnik:Miller}, see also \cite{Moore:PFAsurvey}); hence there cannot be a small basis for the uncountable LOTS in models of CH. 
 However, Moore has shown \cite{5-elements} that assuming PFA there is a five-element basis for the uncountable linear orders: every uncountable linear order will embed one of these five elements.  

Any linear order that embeds into an ordinal must itself be an ordinal, so by a simple argument $\omega_1$ and $\omega_1^*$ will always be minimal order types and so must be included in the basis. Similarly, anything that embeds into a separable linear order must itself be separable (and therefore isomorphic to a set of reals), hence such a linear order must also be in the basis. As we have seen in Section \ref{sec:dense}, Baumgartner proved that under PFA there is a unique (up to isomorphism) $\aleph_1$-dense set of reals, which forms a one-element basis for the uncountable separable linear orders.

Moore developed the theory of those linear orders of size $\aleph_1$ that do not embed $\omega_1$ or $\omega_{1}^{*}$ and have no uncountable separable suborders. These are called Aronszajn lines; they can be constructed as linearisations of Aronszajn trees (uncountable trees with no uncountable branches or levels). Moore proved that under PFA every Aronszajn line contains a Countryman suborder, defined as follows:

\begin{definition}

An {\em Aronszajn line} is a linear order of size $\aleph_1$ that does not contain a suborder of order type $\omega_1$ or $\omega_{1}^{*}$ and has no uncountable separable suborders.
A {\em Countryman line}, $C$, is a linear order of size $\aleph_1$ such that the product $C \times C$ is the union of countably many chains (in the product order).

\end{definition}

The existence of Aronszajn lines can be established in ZFC (see \cite[II]{kunen}). The notion of a Countryman line was first introduced by Countryman in an unpublished article from 1970. Shelah proved that they exist in ZFC \cite{Shelah:countrymanlineproof}. 

\begin{lemma}

Every Countryman line is also an Aronszajn line. 

\end{lemma}
 
$\aleph_1$-dense Countryman lines with a particular property, called non-stationarity, are unique up to isomorphism / reverse isomorphism under PFA.

\begin{definition}

An Aronszajn line, $A$, is {\em non-stationary} if there is a continuous increasing chain of countable subsets of $\omega_1$, $\langle C_\delta : \delta < \omega_1 \rangle$ with union $\omega_1$ such that if $\omega_1$ is the underlying set of $A$ (we assume without loss of generality that it is) then no maximal convex subset of $A \setminus C_\delta$ has end-points.

\end{definition}

For the rest of this section let $X$ and $C$ be, respectively, a fixed $\aleph_1$-dense set of reals and a fixed $\aleph_1$-dense non-stationary Countryman line, with $C^*$ its reverse ordering.

\begin{lemma}
\label{countrymaniso}
\cite[2.1.12]{Todor:WalksonOrdinals}.
Assume PFA, let $C$ and $D$ be $\aleph_1$-dense non-stationary Countryman lines. Then either $C \cong D$ or $C \cong D^*$.

\end{lemma}

Every Countryman line has a non-stationary $\aleph_1$-dense suborder, and can also be extended to a non-stationary $\aleph_1$-dense Countryman line (see \cite{Moore:UniAline}).

\begin{theorem}[Moore, \cite{5-elements}]
\label{Justinbasis}
Assuming PFA, every Aronszajn line contains a Countryman suborder. Consequently, by Lemma \ref{countrymaniso} and the remark immediately following it, $\{X$, $C$, $C^*$, $\omega_1$, $\omega_{1}^{*}\}$ forms a basis for the uncountable linear orders.
Any uncountable linear order must contain an uncountable suborder isomorphic to one of these five.

\end{theorem}

As noted above, Countryman lines can be constructed in ZFC, so by the following result it is clear that this five element basis is in fact minimal for models of ZFC:

\begin{lemma}
Let $C$ be a Countryman line. Then if $D$ embeds into both $C$ and $C^*$, $D$ must be countable.
\end{lemma}

Hence the basis must include both $C$ and $C^*$, or an uncountable suborder of each of them.

However, the five element basis of Theorem \ref{Justinbasis} cannot be a basis for the uncountable LOTS. To ensure that every uncountable linear order {\em continuously} embeds one of the basis elements we need to expand the basis to negate all possible counterexamples to continuity; we will in fact prove that there is an {\em eleven} element basis for the uncountable LOTS. The proof will be divided into three lemmas (Lemmas \ref{Aline-basis}, \ref{sep-basis}, \ref{omega1-basis}) which will establish the additional basis elements needed to ensure that a continuous embedding can always be found. A further lemma (Lemma \ref{min-basis}) proves that this eleven element basis is in fact minimal. 

We will also make use of the following lemma:

\begin{lemma}[Moore, \cite{Moore:UniAline}]
\label{justinlittletheorem}
(PFA) Let $A$ be a non-stationary Aronszajn line, such that $C^*$ does not embed into $A$. Then $A \cong C$.
\end{lemma}

\begin{lemma}\label{Aline-basis}
(PFA) Given an Aronszajn line $A$, one of either $C \times \mathbb{Z}$ or $C^* \times \mathbb{Z}$ must embed into $A$.
\end{lemma}

\begin{proof}
By Lemma \ref{countrymaniso} and Theorem \ref{Justinbasis}, $A$ must embed either $C$ or $C^*$. In the former case, $A$ will embed $C \times \mathbb{Z}$. To see this, note that the $\aleph_1$-density of $C$ implies that $\mathbb{Z}$ embeds into $C$, so $C \times \mathbb{Z} \embeds C \times C$. By Lemma \ref{justinlittletheorem}, $C \times C \cong C$, as none of $C^*, \omega_1, \omega_{1}^{*}$ or $X$ can embed into $C\times C$ and hence it is Aronszajn. Thus $C \times \mathbb{Z}$ embeds into $A$.

If $C^*$ embeds into $A$ then an identical argument tells us that $C^* \times\mathbb{Z}$ also embeds into $A$.
%
%
\end{proof}

Note that the infimum and supremum functions are never well-defined (for infinite decreasing/increasing sequences) on $C\times \mathbb{Z}$ and $C^* \times \mathbb{Z}$, so any embedding from either of these into a given linear order will trivially be continuous. Hence we can replace $C$ and $C^*$ in our basis with the above two orders. However, an analogous result to Lemma \ref{Aline-basis} for the $\aleph_1$-dense set of reals, $X$, strongly fails to hold: if $|B| \geq 2$ then $X \times B \not\embeds X$, for any linear order $B$.
We use instead the following lemma:

\begin{lemma}
\label{sep-basis}

Let $A$ be an uncountable linear order such that there exists a linear order embedding $f : X \rightarrow A$. Then we can find a LOTS embedding $f' : X \times B \rightarrow A$ for some $B \in \{ 1, 2, \omega, \omega^*, \mathbb{Z} \}$.

\end{lemma}

\begin{proof}
First note that under PFA, if $X' \subseteq X$ is $\aleph_1$-dense then $X'$ is isomorphic to $X$.

Let $A$ and $f$ be as in the hypothesis of the lemma. If there is a continuous order-preserving map from $X$ into $A$ (i.e. from $X \times 1$ into $A$) then we are done, so we assume that this is not the case.

In particular then $f$ cannot be continuous, and because every open interval of $X$ is isomorphic to $X$, $f$ cannot be continuous on any convex subset of $X$. This establishes that those points in $X$ where $f$ is discontinuous are dense within $X$; we will show that they are in fact $\aleph_1$-dense in $X$.

Abusing notation slightly, we will write $\mathrm{sup}(f[l(x)]) \not\in A$ to mean the image of the set $l(x)$ has no supremum in $A$, and will likewise write $\mathrm{inf}(f[r(x)]) \not\in A$ when there is no infimum to the image of $r(x)$ in $A$. We define $D$ below as the set of all points in $X$ where $f$ fails to be continuous:

\begin{quote}
 $D = \{ x \in X : \mathrm{sup}(f[l(x)]) \not\in A \text{ or } \mathrm{sup}(f[l(x)]) \neq f(x) \}$
$\text{ }\cup\text{ } \{ x \in X : \mathrm{inf}(f[r(x)]) \not\in A \text{ or } \mathrm{inf}(f[r(x)]) \neq f(x) \}$.
\end{quote}

\begin{claim}

$D$ is $\aleph_1$-dense in $X$.

\end{claim}

\begin{proof}
Assume not. Let $a, b \in X$ be such that $| D \cap (a, b) | \leq \omega$. Then the interval $(a, b)$ is isomorphic to $X$ and $(a, b) \setminus D$ is $\aleph_1$-dense in $(a, b)$. Hence $(a, b) \setminus D$ is isomorphic to $X$ and is continuously embedded into $(f(a), f(b))$ by $f$ (this is because any point where continuity fails for $f \upharpoonright ((a, b) \setminus D)$ would also have to be in $D$, by the density of $(a, b) \setminus D$ within $(a, b)$), which contradicts our assumption.
\end{proof}

\textbf{Continuation of the proof of Lemma \ref{sep-basis}.} We can now classify points in $D$ into four types:

\begin{enumerate}

\item[(i)] Let $D^2 =
\{ x \in D : \mathrm{sup}(f[l(x)]) \in A \text{ and } \mathrm{inf}(f[r(x)]) \in A \}$.\\

\item[(ii)] Let $D^\omega =
 \{ x \in D : \mathrm{sup}(f[l(x)]) \in A \text{ but } \mathrm{inf}(f[r(x)]) \not\in A \}$.\\

\item[(iii)] Let $D^{\omega^*} =
 \{ x \in D :\mathrm{sup}(f[l(x)]) \not\in A \text{ but } \mathrm{inf}(f[r(x)]) \in A \}$.\\

\item[(iv)] Let $D^{\mathbb{Z}} =
 \{ x \in D : \mathrm{sup}(f[l(x)]) \not\in A \text{ and } \mathrm{inf}(f[r(x)]) \not\in A \}$.\\

\end{enumerate}

Clearly (i) - (iv) exhaust all the possibilities for points in $D$, and hence:
\[D = D^2 \cup D^\omega \cup D^{\omega^*} \cup D^\mathbb{Z}\]

Using this we can infer that one of $D^2, D^\omega, D^{\omega^*}, D^\mathbb{Z}$ must be $\aleph_1$-dense in some interval of $X$, given that $D$ itself is. To see this, assume not. Then let $(a, b) \subseteq X$ be such that $| D^2 \cap (a, b) | \leq \omega$. By assumption, $D^\omega$ is not $\aleph_1$-dense in any interval of $X$, so there is an interval $(c, d) \subseteq (a, b)$ such that $| D^\omega \cap (c, d)| \leq \omega$. Repeating this argument two further times gives us a $(g, h) \subseteq X$ such that each of $D^2, D^\omega, D^{\omega^*}, D^\mathbb{Z}$ has countable intersection with $(g, h)$, whereas their union, $D$, has uncountable intersection with $(g, h)$ by $\aleph_1$-density. This is a contradiction.

So one of these four sets must be isomorphic to $X$.
Thus, to finish the proof of Lemma \ref{sep-basis} we must split our argument into four cases.

\textbf{Case 1:} Assume $D^2$ is $\aleph_1$-dense in some interval $(a, b) \subseteq X$, and assume without loss of generality that $D^2 \subseteq (a, b)$. Then $D^2 \cong X$.

The restriction of $f$ to $D^2$ is an injective order-preserving function from $D^2$ into $A$, and we will use it to define an injective continuous order-preserving function $f^\prime$ from $D^2 \times 2$ into $A$. Then composition of $f'$ with some isomorphism $j : X \times 2 \rightarrow D^2 \times 2$ will give a continuous order-preserving injective map from $X \times 2$ into $A$.

To define $f^\prime$, let $x \in D^2$. Then by the definition of $D^2$ we can find a pair of elements $s_x = \mathrm{sup}(f[l(x)])$ and $i_x = \mathrm{inf}(f[r(x)])$, and we set $f^\prime((x, 0)) = s_x$ and $f^\prime((x, 1)) = i_x$. We can do this for all $x \in D^2$; we now need to check that $f^\prime$ is order-preserving and continuous.

Let $(x, y) < (v, w) \in D^2 \times 2$. If $x = v$ then we must have $y = 0$ and $w = 1$. In this case we need only show that 
$f^\prime((x, 0)) < f^\prime((x, 1))$, i.e. $s_x < i_x$. By the definitions of $s_x$ and $i_x$ this is trivial. If however $x < v$ then we need only show that $f^\prime((x, 1)) < f^\prime((v, 0))$. So we need to check that the infimum of $f[r(x)]$ is less than the supremum of $f[l(v)]$ in $A$; this is the case if there is some $c \in X$ contained in both $r(x)$ and $l(v)$ (because then $\mathrm{inf}(f[r(x)]) < f(c) < \mathrm{sup}(f[l(v)])$), which is clearly true by $\aleph_1$-density.

To check that continuity is satisfied by $f^\prime$, let $(x, 0)$ be in $D^2 \times 2$. We need to check that $\mathrm{sup}(f^\prime[l(x) \times 2]) = f^\prime((x, 0))$. But this is clear by the definition of $f^\prime$. Similarly, $f((x, 1))$ is the infimum of $f'[r(x) \times 2]$.

So in this case, $X \times 2$ continuously embeds into $A$.

\textbf{Case 2:}
Assume $D^\omega \cong X$. As before, we use the restriction of $f$ to $D^\omega$ to define a continuous order-preserving map $f^\prime$ from $D^\omega \times \omega$ into $A$. Composition with an appropriate isomorphism gives a continuous order-preserving injection from $X \times \omega$ into $A$.

To define $f^\prime$, let $x \in D^\omega$. Then by the definition of $D^\omega$ we can find an increasing sequence of elements of $A$, $\{ t_{n}^{x} : n < \omega \}$, such that $t_{0}^{x} = \mathrm{sup}(f[l(x)])$ and for all $n < \omega$ we have that for all $y \in r(x)$, $t_{n}^{x} < f(y)$. If we cannot do this then $(f[r(x)])$ must have an infimum in $A$, which contradicts our definition of $D^\omega$. Then we define $f^\prime : D^\omega \times \omega \rightarrow A$ by setting $f^\prime((x, m)) = t_{m}^{x}$.

Checking that $f^\prime$ is order-preserving, injective and continuous is much the same as in Case 1.

\textbf{Case 3:}
 Assume $D^{\omega^*} \cong X$. As would be expected, our definition of a function $f^\prime : D^{\omega^*} \times \omega^* \rightarrow A$ such that $f^\prime$ is order-preserving, injective and continuous proceeds in much the same fashion as in Case 2, except that here we choose $\{ t_{n}^{x} : n < \omega \}$ as a decreasing sequence in $A$, with $t_{0}^{x} = \text{inf}(f[r(x)])$ and for all $n < \omega$ and $y \in f[l(x)]$, $y < t_{n}^{x}$ in $A$. Identify $\omega^*$ with the negative integers and set $f^\prime((x, -m)) = t_{m}^{x}$.

\textbf{Case 4:} 
 If $D^{\mathbb{Z}} \cong X$, we can construct a continuous order-preserving $f^\prime : D^{\mathbb{Z}} \times \mathbb{Z} \rightarrow A$ by taking each $x \in D^\mathbb{Z}$ and choosing a set with order type equal to that of the integers in the interval: 
\[
(\mathrm{sup}(f[l(x)]), \mathrm{inf}(f[r(x)])),
\] which is possible by the definition of $D^\mathbb{Z}$, and mapping $\{ x \} \times \mathbb{Z}$ onto this set in the obvious way. Checking that this gives an $f'$ with the desired properties is much the same as in the previous cases. There is no need to check continuity in this case because no infinite sequence of points in $X \times \mathbb{Z}$ has a limit.

Our construction of $f^\prime$ in all four cases completes the proof of Lemma \ref{sep-basis}.
\end{proof}

To complete our basis for the uncountable LOTS we now need to address those linear orders that only embed the linear order basis elements $\omega_1$ or $\omega_1^*$.

\begin{lemma}
\label{omega1-basis}

Let $\omega_1$ $\bar{\times}$ $\omega^*$ denote the lexicographical product $\prod_{\alpha \in \omega_1} L_\alpha$ where $L_\alpha = 1$ if $\alpha$ is a successor ordinal, and $L_\alpha = \omega^*$ if $\alpha$ is a limit ordinal.
Then if $\omega_1$ embeds into some uncountable linear order $A$, either $\omega_1$ or $\omega_1$ $\bar{\times}$ $\omega^*$ will continuously embed into $A$. Similarly, if $\omega_{1}^{*}$ embeds into $A$ then either $\omega_{1}^{*}$ or $\omega_{1}^{*}$ $\bar{\times}$ $\omega$ (defined analogously) will continuously embed into $A$.

\end{lemma}

\begin{proof}
We prove the first of these two statements; the proof of the second involves only minor modifications of the first.

As before, we assume that there is no continuous map from $\omega_1$ into $A$, but that there is an order-preserving injective function $f : \omega_1 \rightarrow A$.

Let $D = \{ \gamma \in Lim(\omega_1) : f(\gamma) \neq \mathrm{sup}(f[l(\gamma)]) \}$; $D$ must be uncountable, otherwise a final section of $\omega_1$ (equivalently, $\omega_1$ itself) can be mapped continuously into $A$, a contradiction. We must also be able to find an uncountable $D^\prime \subseteq D$ such that $D^\prime = \{ \gamma \in D : \mathrm{sup}(f[l(\gamma)]) \not\in A \}$, because otherwise we could find an $\alpha < \omega_1$ and a continuous map $f^1 : (D \setminus \alpha) \times 2 \rightarrow A$, which of course gives a continuous map from $\omega_1$ into $A$ by the isomorphism of $\omega_1$ and $\omega_1 \times 2$, by setting $f^1((\gamma, 0))$ = $\mathrm{sup}(f[l(\gamma)])$ and $f^1((\gamma, 1)) = f(\gamma)$. This contradicts our assumption that $\omega_1$ does not continuously embed into $A$.

As $D'$ is isomorphic to $\omega_1$, we define a map $f^\prime : D^\prime$ $\bar{\times}$ $\omega^* \rightarrow A$ and show that it is a LOTS embedding, completing the proof. Let $\langle \alpha_i : i < \omega_1 \rangle$ be the increasing enumeration of $D'$, and let $\gamma \in D^\prime$ be such that $\gamma = \alpha_i$ for some limit ordinal $i$. Then we can find a decreasing sequence in $A$, $\{ t_{n}^{\gamma} : n < \omega \}$, with $t_{0}^{\gamma} = f(\gamma)$ and such that for all $m \in \omega$ and $y \in f[l(\gamma)]$ we have $y < t_{m}^{\gamma}$ in $A$. This is clearly possible by the definition of $D^\prime$. 

So we set $f^\prime ((\gamma, -m)) = t_{m}^{\gamma}$ whenever $\gamma = \alpha_i$ for some limit ordinal $i$, where again we are identifying $\omega^*$ with the negative integers, and $f' ((\gamma, 1)) = f(\gamma)$ otherwise. It is easy to check that this map is order preserving and continuity causes no problems since $\omega_1$ $\bar{\times}$ $\omega^*$ has no limit points.
\end{proof}

We have thus proved the following theorem:

\begin{theorem} \label{eleven-elements}

The set:

\begin{align*}
\{ & X, X \times 2, X \times \omega, X \times \omega^*, X \times \mathbb{Z}, \\
& C \times \mathbb{Z}, C^* \times \mathbb{Z}, \\
& \omega_1, \omega_1 \bar{\times} \omega^*, \\
& \omega_{1}^{*}, \omega_{1}^{*} \bar{\times} \omega \}
\end{align*}

forms an eleven element LOTS basis for the uncountable linear orders under PFA.
\end{theorem}

\begin{proof}
To see this, let $A$ be an uncountable linear order. If it embeds $C$ or $C^*$ then it embeds one of $C \times \mathbb{Z}$ or $C^* \times \mathbb{Z}$; trivially, the embedding will be continuous. If it embeds $X$ then Lemma \ref{sep-basis} applies, and if it embeds $\omega_1$ or $\omega_{1}^{*}$ then Lemma \ref{omega1-basis} applies. So it must continuously embed one of the above 11 elements.
\end{proof}

It is possible to establish that this number cannot be improved upon:

\begin{lemma}\label{min-basis}

The 11 element basis is the smallest possible. 
\end{lemma}

\begin{proof}
Well-orders, separable linear orders and Aronszajn lines all exist in ZFC, so the basis must include all of them. It is routine to check that no basis element continuously embeds into any other and there is no uncountable linear order that continuously embeds into two distinct basis elements.
\end{proof}

So we have established that there is a minimal 11 element basis (which cannot be improved upon in any model of ZFC) for the uncountable LOTS.










The existence of this finite basis relies on the inclusion of several linear orders in which no increasing or decreasing sequence has a limit. This is necessary: if not, there will be a bound $\lambda$ on the maximal cofinality of any increasing sequence in any of the basis elements that has a supremum. But there are linear orders such that every point has cofinality $\kappa$, for any $\kappa$. So when $\lambda < \kappa$ such a linear order could not continuously embed any basis element that has well-defined suprema. Similarly for decreasing sequences.

Each point in a linear order is both an infimum and a supremum (of a decreasing/increasing sequence) if and only if the linear order is dense. Can there be a small basis for dense linear orders? By the argument in the above paragraph we have to restrict attention to linear orders where every point has small cofinality/coinitiality, again under PFA.

\begin{theorem}
(PFA) There is a six element basis for those uncountable dense LOTS in which all points have cofinality and coinitiality $\omega$. The basis is $\{ X, X \times \mathbb{Q}, C \times \mathbb{Q}, C^* \times \mathbb{Q}, \omega_1 \times \mathbb{Q}, \omega_{1}^{*} \times \mathbb{Q} \}$.
\end{theorem}

\begin{proof}
We begin by proving the following:

\begin{claim}
If $L$ is dense and all points in $L$ have cofinality and coinitiality $\omega$ then $\mathbb{Q}$ embeds continuously into $L$.
\end{claim}
\begin{proof}
We find a countable dense subset $L' \subseteq L$ without endpoints such that for $x \in L'$, sup$(l(x) \cap L')$ $=x=$ inf$(r(x) \cap L')$. Then $L' \cong \mathbb{Q}$ and the latter condition guarantees continuity.
We can find such an $L'$ by simply iterating a process of choosing a point and then choosing countable increasing and decreasing sequences which have this point as their limit. Continue this until the order is dense and without endpoints, in which case it must be isomorphic to $\mathbb{Q}$.
\end{proof}

To prove the Theorem, we argue as in Lemmas \ref{Aline-basis}, \ref{sep-basis}, \ref{omega1-basis}. If $A$ is a dense uncountable linear order it must embed one of the five basis elements. Again, we assume these embeddings are not continuous and find a set of points of discontinuity that are isomorphic to an element of the basis. By the density of $L$ we can continuously embed a copy of $\mathbb{Q}$ into the interval $(f[l(x)], f[r(x)])$, where $x$ is a point of discontinuity. It is easy to see that this works. Note that density requires us to include the full lexicographical product $\omega_1 \times \mathbb{Q}$ rather than $\omega_1$ $\bar{\times}$ $\mathbb{Q}$ as in the previous theorem, and similarly for $\omega_{1}^{*}$.

\end{proof}

\section{Open questions}\label{sec:questions}

We collect here some open questions arising from this paper:

\begin{question}
Shelah, in \cite{sh100}, has proved it consistent that there is a universal linear order at $\aleph_1$ for the case where $\aleph_1 < 2^{\aleph_0}$ holds. However, this is not a universal LOTS. Is it consistent that there is a universal LOTS in this situation?

\end{question}

\begin{question}
Can there ever be a universal LOTS for uncountable $\kappa$?
\end{question}

\begin{question} In \cite{Moore:UniAline}, Moore proves that there is a universal Aronszajn line assuming PFA. This is not a universal LOTS. Can it be shown that there is no universal Aronszajn LOTS in any model of ZFC?
\end{question}

\begin{question}
Can there be a finite (or even countable) basis for the $\aleph_1$-dense LOTS in which all points have countable cofinality and coinitiality?
\end{question}

\bibliographystyle{amsplain}

\providecommand{\bysame}{\leavevmode\hbox to3em{\hrulefill}\thinspace}
\providecommand{\MR}{\relax\ifhmode\unskip\space\fi MR }
\providecommand{\MRhref}[2]{%
  \href{http://www.ams.org/mathscinet-getitem?mr=#1}{#2}
}
\providecommand{\href}[2]{#2}

\end{document}